\documentclass[11pt]{amsart}

\newtheorem{thm}{Theorem}[section]
\newtheorem{lm}[thm]{Lemma}

\theoremstyle{definition}

\newtheorem*{df*}{Definition}

\theoremstyle{remark}

\newtheorem*{rem*}{Remark}

\numberwithin{equation}{section}

\newcommand{\ci}[1]{_{ {}_{\scriptstyle #1}}}

\newtheorem{zamech}{Remark}

\newcommand{\cB}{\mathcal{B}}
\newcommand{\cD}{\mathcal{D}}

\newcommand{\cL}{\mathcal{L}}

\newcommand{\f}{\varphi}

\newcommand{\R}{\mathbb{R}}

\newcommand{\p}{\partial}

\newcommand{\wt}{\widetilde}

\newcommand{\al}{\alpha}

\newcommand{\La}{\langle }
\newcommand{\Ra}{\rangle }

\newcommand{\sd}{{\scriptstyle\Delta}}
\newcommand{\bu}{\mathbf{u}}
\newcommand{\bn}{\mathbf{n}}

\newcommand{\bw}{\mathbf{w}}
\newcommand{\bff}{\mathbf{f}}

\newcommand{\bT}{\mathbf{T}}

\def\cyr{\fontencoding{OT2}\fontfamily{wncyr}\selectfont}
\DeclareTextFontCommand{\textcyr}{\cyr}

\begin{document}

\title[Carleson--Buckley measures]{Carleson--Buckley measures beyond the scope of $A_\infty$ and their applications}
\author{F. Nazarov, A. Reznikov, S. Treil, and A. Volberg}

\maketitle

\section{Introduction}
\label{intro}

Carleson measures are ubiquitous in Harmonic Analysis.  In the paper of Fefferman--Kenig--Pipher \cite{FKP} an interesting class of Carleson measures was introduced for the need of regularity problems of elliptic PDE. These Carleson measures were associated with $A_\infty$ weights. In discrete setting (we need exactly discrete setting here) they were studied in  Buckley's \cite{Buc1}, where they were associated with dyadic $A_\infty^d$.

Our goal here is to show that such Carleson--Buckley measures (in discrete setting) exists {\it for virtually any} positive function (weight). Of course some modification is needed, because it is known  (see below) that Carleson property of Buckley's measure are {\it equivalent} to the weight to be in $A_\infty^d$.

However a very natural generalization of those facts exist for {\it any} weight, and of course, as a natural application to special case $w\in A_\infty^d$ it 
gives Buckley's results. The same can be said for continuous version of \cite{FKP}.

Notice that these Carleson--Buckley measures for {\it general} weight immediately gave some applications. In preprints \cite{AL}, \cite{NRTV1}, \cite{NRTV2} (\cite{NRTV2} is a full text version of a short \cite{NRTV1}) a rather long standing problem called ``bump condition problem" has been solved. The methods are different, and \cite{NRTV1}, \cite{NRTV2}  formulated only in metric $L^2$, but \cite{NRTV1}, \cite{NRTV2}  are based on these general Carleson--Buckley measure pertinent to {\it general} (not $A_\infty$) weight. In this sense \cite{NRTV1}, \cite{NRTV2}  are slightly more general than \cite{AL}. Although it is feasible that our Carleson--Buckley results of \cite{NRTV1}, \cite{NRTV2}  and in this note also can be obtained by ``local mean oscillation decomposition" of \cite{AL}.

\subsection{$A_\infty$ and Carleson measures associated with it}
\label{A1}

We formulate two results, the first belongs to Buckley, the second is probably a folklore one.
Let $w\in A_\infty^d$, let $I\in \cD$ be a dyadic interval, $I_+, I_-$ being its right and left halves, 
$\Delta_I w:= \La w\Ra_{I_+}-\La w\Ra_{I_-}$,  then

\begin{equation}
\label{Buc}
\forall J\in \cD\,,\,\,\sum_{I\subset J, \,I\in \cD}\frac{|\Delta_I w|^2}{\La w\Ra_I} |I| \le C w(J)\,,
\end{equation}
where $C$ depends only on $A_\infty$ characteristic  $[w]_{A_\infty^d}$ of $w$. It is known how exactly it depends on $A_\infty$ characteristic  $[w]_{A_\infty^d}$ of $w$, see \cite{BR}. But we do not need this here.

Clearly this is the statement about $w$-Carleson measure associated with $A_\infty^d$ weight. Another statement 
is the following: again  let $w\in A_\infty^d$, $\{\al_I\}_{I\in \cD}$  be a dyadic Carleson sequence (discrete measure), meaning that $\al_I\ge 0$ and 

\begin{equation}
\label{C1}
\forall J\in \cD\,,\,\,\sum_{I\subset J, \,I\in \cD}\al_I |I| \le C|J|\,,
\end{equation}
where the best $C$ is called the Carleson norm of the sequence (measure) $\{\al_I\}_{I\in\cD}$.
Then
\begin{equation}
\label{Buc2}
\forall J\in \cD\,,\,\,\sum_{I\subset J, \,I\in \cD}\La w\Ra_I \al_I|I|  \le C w(J)\,,
\end{equation}

This second folk result we will even prove now: as $\La w\Ra_I \le C\La w^{1/2}\Ra_I^2$ (it is in $A_\infty^d$ so it satisfies the Reverse H\"older Inequality (RHI)), we can replace \eqref{Buc2} by
$$
\sum_{I\subset J, \,I\in \cD}\La f\Ra_I ^2\al_I|I|  \le C\int_J f^2 dx\,,
$$
where $f:= w^{1/2}$. But the latter inequality is well known (see Garnett's book, for example) property of Carleson sequences. We are done with \eqref{C1}. See \cite{Buc1} and \cite{FKP} for \eqref{Buc}.

\section{Not $A_\infty$ weights}
\label{not}

What if $w$ is not an $A_\infty$ weight?
Obviously \eqref{Buc}, \eqref{Buc2} fail. It is known that \eqref{Buc} implies $A_\infty^d$, and we  even know the sharp dependence of $[w]_{A_\infty^d}$ on $C$ from \eqref{Buc}.  In what concerns \eqref{Buc2}, the maximal operator is not bounded in $L^1$ unfortunately.

But we are going to prove the following theorems for general $w$. They will of course imply the previous section trivially. But they also implied the bump conjecture, see \cite{NRTV1}, \cite{NRTV2} (we already mentioned another solution in \cite{AL}). See also immense amount of references in these papers and in the book \cite{CU-M-P-book}. Moreover, the theorems below have stronger versions discussed in the last Section. Therefore, these theorems prove the bump conjecture under weaker ``bump" assumptions.

To formulate the results we need Orlicz norms, and actually, something else expressed by function $\Psi$ below.

\subsection{Orlicz norms and distribution functions}
\label{orlicz}

\subsection{A lower bound for the Orlicz norm}
Let $\Phi$ be a continuous non-negative increasing convex function such that $\Phi(0)=0$ and 
$\int^{+\infty}\frac{dt}{\Phi(t)}<+\infty$. Define $\Psi(s)$ parametrically by $\Psi(s)=\Phi'(t)$ when 
$s=\frac{1}{\Phi(t)\Phi'(t)}$ ($t>0$). Then $\Psi(s)$ is positive and decreasing for 
$s>0$ and $s\Psi(s)$ is increasing. Moreover $\int_0\frac{ds}{s\Psi(s)}<+\infty$. Indeed, using our parameterization we can rewrite the last integral as 
$$
\int^{+\infty}\left(\frac{1}{\Phi(t)}+\frac{\Phi''(t)}{\Phi'(t)^2}\right)\,dt\,.
$$
The first integral converges by our assumption and the second integrand has a bounded near $+\infty$
antiderivative $\frac{-1}{\Phi'(t)}$.

Let $w\ge 0$ on  $I\subset X$. Define the normalized distribution function $N$ of $w$ by
$$
N(t)=N_I^w(t)=\frac{1}{\mu(I)} \mu(\{x\in I:w(x)>t\})
$$

\begin{lm}
\label{l:Orl_LB}
Let $\Psi : (0,1]\to \R_+$ be a decreasing function such that the function $s\mapsto s\Psi(s)$ is increasing. Let $\Phi$ be a Young function and let
\begin{align*}
\Psi(s) \le C \Phi'(t) \qquad \text{where} \quad s = \frac{1}{\Phi(t) \Phi'(t)}
\end{align*}
for all sufficiently large $t$. Then for $N=N_I^w$
\begin{align}
\label{n(N)}
\bn\ci\Psi (N) := \int_0^\infty N(t)\Psi(N(t))\,dt\le C\|w\|\ci{L^\Phi(I)}\,.
\end{align}
\end{lm}

\begin{proof}
The left hand side scales like a norm under multiplication by constants, so it is enough to 
show that if $\|w\|_{L_\Phi(I)}\le 1$, i.e.,
\[
\frac{1}{|I|}\int_I\Phi(w)=\int_0^\infty N(t)\Phi'(t) \,dt \le 1
\]
then $\bn\ci\Psi(N)$ is bounded by a constant. Since $s\Psi(s)$ increases, we may have 
trouble only at $+\infty$ It is cleat that it suffices to estimate the
integral over the set
where $\Psi(N(t))>\Phi'(t)$ but since $\Psi$ is decreasing this means that $N(t)\le C/{ (\Phi(t)\Phi'(t)) }$, so we get
at most $\int^{+\infty}\Phi(t)^{-1} dt$ and we are done.
\end{proof}

\begin{rem*}
In the above Lemma \ref{l:Orl_LB} we do not need the assumption that 
\begin{align}
\label{integrability_01}
\int_0 \frac{1}{s\Psi(s)} ds <\infty.  
\end{align}
But in what follows this assumption will be needed, and the reasoning in the beginning of this section shows that for any Young function $\Phi$ satisfying $\int^\infty (\Phi(t))^{-1} dt <\infty$ we can find $\Psi$ from Lemma \ref{l:Orl_LB} satisfying \eqref{integrability_01}.  
\end{rem*}

\subsection{Examples} In the above section only the behavior of $\Phi$ at $+\infty$ (equivalently, the behavior of $\Psi $ near $0$) was important, so we will concentrate our attention there. 

Let $\Phi(t) = t (\ln t)^\alpha$, $\alpha>1$ near $\infty$. Then 
\[
\Phi'(t) \sim (\ln t)^{\alpha}, \qquad \Phi(t)\Phi'(t) \sim t (\ln t)^{2\alpha} , 
\]
so $\Psi(s) := (\ln (1/s))^\alpha$ satisfies the assumptions of Lemma \ref{l:Orl_LB}: to see that we  notice
\[
\ln(\Phi(t)\Phi'(t) )\sim \ln t. 
\]

If $\Phi(t) = t \ln t (\ln\ln t)^\alpha$, $\alpha>1$, then
\[
\Phi'(t)\sim  \ln t (\ln\ln t)^\alpha, \qquad \Phi(t)\Phi'(t) \sim t (\ln t)^2 (\ln\ln t)^{2\alpha}
\]
and $\Psi(s) = \ln(1/s) (\ln\ln (1/s))^\alpha$ works. because again $\ln(\Phi(t)\Phi'(t) )\sim \ln t$. 

Note that in both examples $\int_0 (s\Psi(s))^{-1} ds <\infty$. 

The examples of Young functions with higher order logarithms are treated similarly.

\subsection{Differential Embedding Theorems with weight not satisfying $A_\infty$}
\label{DE}

Here is the analog of Buckley's inequality \eqref{Buc} for weights without $A_\infty$ property.
\begin{thm}
\label{BUC}
Let $\Psi$ be as in the previous subsection. Then for any weight $w$  such that $\bn\ci \Psi(N_I^w)<\infty$ for all $I\in \cL$
\begin{align}
\label{d-embed}
\sum_{I\in\cD, I \subset J} \bn\ci \Psi(N_I^w)^{-1} (\Delta_{I}w)^2 \le C w(J) \,,
\end{align}
 here in the summation we skip $I$ on which $w\equiv0$. 
%
\end{thm}

It is well known that the previous theorem will imply the following {\it differential embedding theorem}, which was instrumental in the solution of the bump conjecture in \cite{NRTV1}, \cite{NRTV2}.

\begin{thm}
\label{t:fd-embed}
Let $\Psi$ be as in the previous subsection. Then for any weight $w$  such that $\bn\ci \Psi(N_I^w)<\infty$ for all $I\in \cL$
\begin{align}
\label{d-embed-01}
\sum_{I\in\cD, I\subset J} \bn\ci \Psi(N_I^w)^{-1} (\Delta_{I}(fw))^2 \le C\int_J f^2 w dx \,,
\end{align}
for all $f\in L^2(w)$; here in the summation we skip $I$ on which $w\equiv0$. 
%
\end{thm}

\subsection{Embedding theorem with weight not satisfying $A_\infty$}
\label{E}

Another theorem, which was instrumental in the solution of the bump conjecture in \cite{NRTV1}, \cite{NRTV2}, is the analog of \eqref{Buc2} for weight not satisfying any $A_\infty$ conditions.

\begin{thm}
\label{C2}
Let $\{\al_I\}_{I\in \cD}$ be a Carleson sequence as in \eqref{C1}. Let $\Psi$ be as in subsection \ref{orlicz}.  Then for any weight $w$  such that $\bn\ci \Psi(N_I^w)<\infty$ for all $I\in \cL$
\begin{align}
\label{embed}
\sum_{I\in\cD, I\subset J} \bn\ci \Psi(N_I^w)^{-1} \La w\Ra_I^2 \al_I\,|I| \le Cw(J) \,,
\end{align}
here in the summation we skip $I$ on which $w\equiv0$. 
%
\end{thm}

It is well known that the previous theorem will imply the following {\it differential embedding theorem}, which was instrumental in the solution of the bump conjecture in \cite{NRTV1}, \cite{NRTV2}.

\begin{thm}
\label{C22}
Let $\{\al_I\}_{I\in \cD}$ be a Carleson sequence as in \eqref{C1}. Let $\Psi$ be as in subsection \ref{orlicz}.  Then for any weight $w$  such that $\bn\ci \Psi(N_I^w)<\infty$ for all $I\in \cL$
\begin{align}
\label{embed2}
\sum_{I\in\cD, I\subset J} \bn\ci \Psi(N_I^w)^{-1} \La fw\Ra_I^2 \al_I\,|I| \le C\int_J f^2 w dx \,,
\end{align}
here in the summation we skip $I$ on which $w\equiv0$. 
%
\end{thm}

\begin{zamech} The reader should notice that this is obtained by (strangely enough) writing in \eqref{Buc2} $\La w\Ra_I$ as 
$\frac{\La w\Ra_I^2}{\La w\Ra_I}$, and replacing the denominator by {\bf bigger} ``bumped" average.
The same happened in Theorem \ref{BUC}: we took Buckley's inequality \eqref{Buc} and replaced its denominator by {\bf bigger} ``bumped" average, by the way, the same one.
\end{zamech}

\section{Proofs of Differential Embedding Theorems}
\label{PDE}

Let  $\phi(s):=s\Psi(s)$. Let $U(s)=\frac 1{\phi(s)}$. Define $B(s)$ on $[0,1]$ by $B(0)=B'(0)=0$, $B''(s)=U(s)$. Since $U$ is integrable, $B$ is well-defined and $B(s)\le Cs$. Since $U$ is decreasing, we have the finite difference 
inequality
$$
\frac 12[B(w+\Delta w)+B(w-\Delta w)]\ge B(w)+\frac 12 U(w)(\Delta w)^2\,.
$$
Consider 
$$
\mathcal B(I)=\int_0^\infty B(N(t))\,dt\le C\int_0^\infty N(t)\,dt
=C\frac{1}{|I|}\int_I w\,.
$$
Note that 
$$
\frac 12[\mathcal B(I_-)+\mathcal B(I_+)]\ge \mathcal B(I)+
\frac 12\int_0^\infty \frac 1{\phi(N(t))}(\Delta N(t))^2\,dt.
$$
By Cauchy-Schwarz, this integral is at least
$$
\Bigl[\int_0^\infty \phi(N(t))\,dt \Bigr]^{-1}
\Bigl[\int_0^\infty \Delta N(t)\,dt \Bigr]^2\,.
$$
But the first integral is dominated by $\|w\|_{L_\Phi(I)}$ and the second one is also known as 
$\Delta_I w$. Thus, we have the Bellman function proof of the fact that the sequence
$\frac{(\Delta_I w)^2}{\|w\|_{L_\Phi(I)}}|I|$ is $w$-Carleson.  So Theorem \ref{BUC} is already proved. We promised to deduce Theorem \ref{t:fd-embed} from it.

\subsection{Proof of Theorem \ref{t:fd-embed}}

Another, more Bellman technique proof is in \cite{NRTV2}. 

Let $w,v$ be any positive weights. Then 
$$
\sum_{I}\frac{|\Delta_I(fw)|\cdot|\Delta_I(gv)|}{\sqrt{\|w\|_{L_\Phi(I)}\|v\|_{L_\Phi(I)}}}|I|
\le C\|f\|_{L^2(w)}\|g\|_{L^2(v)}\,.
$$

Indeed, let us use our lovely shifted Haar functions $h_I^w$ and $h_I^v$ normalized
in $L^2(w)$ and $L^2(v)$ respectively and write, as usual,
$$
\Delta_I(fw)=\alpha_I(f,h_I^w)_{L^2(w)}/\sqrt{|I|}+
\frac{\langle fw \rangle_I}{\langle w \rangle_I}\Delta_I\,w
$$
with $|\alpha_I|\le\sqrt {\langle w \rangle_I}$
and, similarly,
$$
\Delta_I(gv)=\beta_I(f,h_I^v)_{L^2(v)}/\sqrt{|I|}+
\frac{\langle gv \rangle_I}{\langle v \rangle_I}\Delta_I\,v\,.
$$
As usual, we have four sums to estimate. Using the fact that the $L_\Phi$ norm dominates the $L^1$ norm,
we see that the first sum is fine by the standard Parceval inequality. To bound the second sum, we can ignore the factor $\frac{|\beta_I|}{\sqrt{\|v\|_{L_\Phi(I)}}}\le 1$, use Cauchy--Schwarz inequality, and bound  
the sum
$$
\sum_I
\left[\frac{\langle fw \rangle_I}{\langle w \rangle_I}\right]^2\frac{(\Delta_Iw)^2}{\|w\|_{L_\Phi(I)}}|I|
$$
by $\|f\|_{L^2(w)}^2$ using the classical weighted Carleson embedding theorem. The third sum is similar. At last, applying Cauchy--Schwarz inequality to the fourth sum, we see that we again can use the same Carleson type bounds but now for both 
$w$ and $v$.

\section{Proofs of  Embedding Theorems}
\label{PE}

We can think that Carleson constant in \eqref{C1} is $1$. Let
$$
A:= A_I :=\frac{1}{|I|} \sum_{I'\subset I, I'\in \cD} \al_I |I| \le 1\,.
$$
So in what follows we can think that $0\le A\le 1$.

Consider 

$$
T(A, N) := N\int_0^{N/(A+1)} \frac1{\phi(s)}\,ds\,.
$$

Here we at least used that $1/\phi(s)$ is integrable at $0$. Notice that $0\le T(A, N) \le cN$, and that
\begin{equation}
\label{der}
-\frac{\partial T}{\partial A} \ge c\frac{N^2}{\phi(N)}
\end{equation}
because of the doubling condition on $\phi$ and because $0\le A\le 1$.


We need to check that $T$ is {\bf convex}.  

\vspace{.1in}

\noindent {\bf Lemma.}
\label{always}
Function $T$ is convex  for any $\phi$ because
\begin{equation}
\label{usl}
\frac{\phi'(s)}{\phi(s)} \le \frac{2}{s}\,,\, \forall s\in (0,1)\,.
\end{equation}
In fact,
 for our function even a stronger inequality $s\phi'(s) \le \phi(s)$ is satisfied.  Indeed, since $\phi(s) = s\Psi(s) $ is increasing and $\Psi$ is decreasing, then 
\begin{align*}
0\le \phi'(s)=(s\Psi(s))' = \Psi(s) + s \Psi'(s) \le \Psi(s) =\frac{\phi(s)}{s}\,,
\end{align*}
(the second inequality holds because $\Psi $ is decreasing).

\begin{proof}
Notice that function
$$
f(x,y) = x G(\frac{x}{y})
$$ always satisfies Monge-Amp\`ere equation: $\partial^2_{xx} f \partial^2_{yy} f -(\partial^2_{xy} f)^2 =0$. In fact, such a function is linear on foliating lines  $y=cx$. One can of course make it by direct calculation as well. So our $T$ is such.

Now let us compute the second derivative of $T$ in $A$. It is 

$$
2\frac{N^2}{\phi(N/(A+1)) (A+1)^3} - \frac{N^3}{(A+1)^4} \frac{\phi'(N/(A+1))}{\phi(N/(A+1))}=
$$
$$
\frac{N^2}{(A+1)^3 \phi(N/(A+1))} \bigg( 2- \frac{N}{A+1} \frac{\phi'(N/(A+1))}{\phi(N/(A+1))}\bigg)\ge 0\,,
$$
because of \eqref{usl}.

 Now because the Hessian's determinant is zero, we get automatically that the second derivative in $N$ is also positive. We are done.

\end{proof}

Now we define $B(A,N)=CN-T(A,N)$, and 
$$
B(I) := \int_0^{\infty}B(A_I, N_{w,I}(t))\,dt\,,
$$
where $A_I= \frac1{I}\sum_{\ell\in \mathcal{D},\, \ell\subset I} \alpha_{\ell}\ell$, and
$N(t):=N_{w, I}(t):= \frac1I\{x\in I: w(x)>t\}$.

By the properties of $B(A, N)$ (concavity and \eqref{der}) the following inequality holds for any $I$ and any $t\in (0, \infty)$:
\begin{equation}
\label{main}
B(I) -\frac12(B(I_+)+B(I_-)) \ge c \alpha_I\int_0^{\infty}\frac{N^2(t)}{\phi(N(t))}\, dt\ge c\alpha_I\frac{(\int N(t)dt)^2}{\int \phi(N(t))dt}\,.
\end{equation}

The latter fraction is at least
$$
c\alpha_I \frac{\langle w\rangle_I^2}{\|w\|_{L_\Phi,I}}\,,
$$
and we can follow the usual steps of Bellman induction to get \eqref{Buc2}, \eqref{embed}.

\section{Application}
\label{appl}

We want to sketch the application from \cite{NRTV2} to the so-called ``bump conjecture".

\subsection{Bellman function and main differential inequality}
\label{s:Bell01-DE-01}
Let $\f(s):={s\Psi(s)}$. Multiplying $\Psi$ by an appropriate constant we can assume without loss of generality that 
\begin{equation}
\label{norm-Psi}
\int_0^1 \frac1{\f(s)} ds =1.
\end{equation} 

 Define $m(s)$ on $[0,1]$ by $m(0)=m'(0)=0$, $m''(s)=1/\f(s)$. Identity \eqref{norm-Psi} implies that $m$ is well-defined and $m'(s) \le 1$,  $m(s)\le s$. 
For a distribution function $N=N_I^w$ define
\begin{align}
\label{u(N)}
\bu (N) = \int_0^\infty (2 N(t) - m(N(t))) dt = 2 \La w\Ra\ci I - \int_0^\infty m(N(t)) dt.
\end{align}

For the scalar variable $\bff\in \R$ and  the distribution function $N$ define the Bellman function $\wt\cB(\bff, N)=\cB(\bff, \bu(N))$ where 
\[
\cB(\bff,\bu ) = \frac{\bff^2}{\bu}. 
\]

\subsection{Main inequality 
in the finite difference form}

It is proved in \cite{NRTV2} that if $\bn:= \int_0^{\infty} \phi (N(t))\, dt$, then by computing second derivative of $\wt\cB$ in the direction $\sd = (\sd \bff, \sd N)$ we get 
\begin{align}
\label{MainIneq2} 
\wt\cB_{\sd}'' \ge \frac{2(\sd\bff)^2}{2\kappa^2\bn + \bu} \ge \frac{2(\sd\bff)^2}{2\cdot 2^2\bn + \bu} \ge c \frac{(\sd\bff)^2}{\bn}\,.
\end{align}
This implies the following.

\begin{lm}
\label{l:MainIneq}
Let 
\begin{align*}
\bff = \frac{\bff_1 + \bff_2}{2}, \qquad N(t) = \frac{N_1(t) + N_2(t)}{2}.
\end{align*}
Then 
\begin{align}
\label{d-MainIneq}
\frac12 \Bigl( \cB(\bff_1, \bu(N_1))   + \cB(\bff_2, \bu(N_2)) \Bigr) - \cB(\bff, \bu(N)) \ge \frac{c}4   \cdot  \frac{ (\bff_1 -\bff)^2}{\bn(N)}. 
\end{align}
for some positive absolute constant $c$. (Note that $\bff_1 -\bff = \bff-\bff_2$, so we can replace $(\bff_1 -\bff)^2$ in the right side by $(\bff_2 -\bff)^2$)
\end{lm}
\begin{proof}
Notice that 
\begin{align}
\label{conc1}
\frac{s_1+s_2}{2} \Psi\left( \frac{s_1+s_2}{2} \right) \ge \frac{s_1+s_2}{2} \Psi\left( s_1+s_2 \right)
\ge\frac12 s_1 \Psi(s_1); 
\end{align}
here  the first inequality holds because $\Psi$ is decreasing and the second one because $s\Psi(s)$ is increasing. Of course, we can interchange $s_1$ and $s_2$ in the above inequality. 

Let $\sd\bff := \bff_1 -\bff$, $\sd N :=N_1-N$. 
Define
\begin{align*}
F(\tau) = \cB(\bff + \tau\sd\bff, \bu(N +\tau\sd N)) + \cB(\bff - \tau\sd\bff, \bu(N -\tau\sd N)) 
\end{align*}
Taylor's formula together with the estimate  (3.8) from \cite{NRTV2} imply that
\begin{align}
\label{d-conv-1}
F(1) - F(0) \ge \frac{c}2 (\sd\bff)^2 \left( \frac1{\bn(N+\tau\sd N)} + \frac1{\bn(N -\tau\sd N)} \right)
\end{align}
for some $\tau\in (0,1)$. 

Estimate \eqref{conc1} implies that 
\[
\bn(N) \ge \frac12 \bn (N\pm\tau\sd N), 
\]
so 
\[
\left( \frac1{\bn(N+\tau\sd N)} + \frac1{\bn(N -\tau\sd N)} \right) \ge \frac1{\bn(N)}. 
\]
Then it follows from \eqref{d-conv-1} that 
\[
F(1) - F(0) \ge \frac{c}2 \cdot \frac{ (\sd\bff)^2}{\bn(N)}. 
\]
Recalling the definition of $F$ and dividing this inequality by $2$ we get \eqref{conc1}.  
\end{proof}

\subsection{General case}

Let $\f$ and $\wt\cB$ be as above. 
\begin{lm}
\label{maindiff2}
Let $\bff, \bff_k \in\R$, $\alpha_k\in \R_+$ and the distribution functions $N$, $N_k$, $k=1, 2, \ldots, n$ satisfy
\[
\bff = \sum_{k=1}^n \alpha_k \bff_k, \qquad N = \sum_{k=1}^n \alpha_k N_k,  \qquad 
\sum_{k=1}^n \alpha_k =1.  \ 
\]
Then 
\[
-\wt\cB(\bff, N) + \sum_{k=1}^n \alpha_k \wt\cB(\bff, N_k) \ge \frac{c}{16}\cdot \frac1{\bn(N)} \left( \sum_{k=1}^n \alpha_k | \bff_k -\bff |  \right)^2 
\]
\end{lm} 

See the explanation how this lemma follows from the previous one in \cite{NRTV2}. Of course it is an interesting exercise in convexity. Essentially, boldface variables should be substituted by averages of $w$ over a dyadic interval $I$ and over its $2^n$ children of $n$-th generation. Distribution function $N$ should be thought as the normalized distribution function of $w$ on $I$, and $N_k$'s are normalized distribution  functions of $w$ on the children of $n$-th generation.

On the other hand,  Lemma \eqref{maindiff2}  with $\al_k=2^{-n}$ {\it exactly} means the boundedness of a dyadic shift of complexity $n$ (actually of any ``slice" of dyadic shift of complexity $n$ as one can see in \cite{NRTV2}), it is, in fact, the required inequality in a different language, see why is that in \cite{NRTV2}, but it is really just a simple observation!

This is how we deal with dyadic shifts, because every shift of complexity $n$ has $n$ slices (see \cite{NRTV2}), the estimate becomes linear in complexity.

Now there are also paraproducts to be treated. Here we state the embedding theorem, which gives the estimate for the paraproduct operator in two-weight situation under the bump condition.

\begin{thm}
\label{t.embed-bump}
Let $\f: (0,1]\to \R_+$ be a bounded  increasing function such that $s\mapsto \f(s)/s$ is decreasing, $\f(s) \ge s$ and 
\begin{align*}
\int_0^1\frac{1}{\f(s)} ds <\infty. 
\end{align*}
For any normalized Carleson sequence $\{a\ci I\}\ci{I\in\cD}$ ($a\ci I\ge 0$), i.e.~for any sequence satisfying  
\[
\sup_{I\in\cD} |I|^{-1}\sum_{I'\in\cD: I'\subset I} a\ci I' |I'| \le 1
\]
we get
\[
\sum_{I\in\cD} \frac{\La fw\Ra\ci I^2}{\bn(N_I^w)} a\ci I \le C \|f\|_{L^2(w)}^2
\]
\end{thm}

\subsection{An auxiliary function}
Let $\f$ be as above in Theorem \ref{t.embed-bump}. 
For the numbers $A\in[1, 2]$, $N\in \R_+$ define 
\[
T(A, N) := N \int_0^{N/A} \frac{1}{\f(s)} ds
\]

This is our function from the previous section that gives Carleson estimate in \eqref{Buc2}.

\subsection{Bellman function and the main differential inequality. }
Let now $N$ be a distribution function, and let 
\[
\bT (A, N) = \int_0^\infty T(A, N(t)) \, dt . 
\]

As in Section \ref{s:Bell01-DE-01} assume, multiplying $\f$ by an appropriate constant,  that
\[
\int_0^1 \frac{1}{\f(s)} ds =1.
\]
Then $T(A, N(t)) \le N(t)$, so
\[
\bT(A, N) \le \int_0^\infty N(t) dt =:\bw =\bw(N). 
\]

For $\bff\in\R$, $M\in[0,1]$ and for a distribution function $N$ define 
$
\wt\cB(\bff, N, M) : = \cB(\bff, \bu(N, M))
$, 
where 
\[
\cB(\bff, \bu) = \frac{\bff^2}{\bu}
\]
and
\begin{align*}
\bu =\bu(N, M) & = 2 \int_0^\infty N(t) dt - \bT(N, M+1)
\\
& =: 2 \bw(N) - \bT(N, M+1). 
\end{align*}
Note that $\bu(N) \ge \bw(N)$. 

Exactly as before, we get 
\begin{align}
\label{dB/dM-02}
-\frac{\p \wt\cB}{\p M} \ge \frac1{16}\cdot \frac{\bff^2}{\bn}
\end{align}
This inequality (together with the convexity of $\wt\cB$) is the main differential inequality for our function. 

\subsection{Finite difference form of the main inequality}
Let $X=(\bff, N, M)$, $X_k =(\bff_k, N_k, M_k)$, ($\bff, \bff_k\in \R$, $M, M_k\in[0,1]$, $N$, $N_k$ are the distribution functions) satisfy 
\begin{align*}
\bff = \sum_{k=1}^n \alpha_k \bff_k, \qquad N = \sum_{k=1}^n \alpha_k N_k, \qquad M = a + \sum_{k=1}^n \alpha_k M_k, \ a\ge 0, 
\end{align*}
where 
\begin{align*}
 \sum_{k=1}^n \alpha_k =1, \qquad \alpha_k \ge 0. 
\end{align*}
Then
\begin{align}
\label{discr-MainIneq2-03}
- \wt \cB(X) + \sum_{k=1}^n \alpha_k \wt\cB(X_k) \ge \frac1{16} \cdot \frac{a\bff^2}{\bn} 
\end{align}
where $\bn=\bn(N)$.

Again, see \cite{NRTV2} for the explanation that this main inequality  is {\it exactly} the boundedness of the paraproduct given a bump condition.

\section{Discussion}
\label{di}

The reader already probably noticed that we do not require
 the bump condition. Instead of $\|\cdot \|_{\Phi, I}$ and the requirement that the product of such quantities is bounded, we require that the product of $n_{\Psi}(\cdot)$ is bounded.  We saw in Lemma \ref{l:Orl_LB} that
 \begin{align}
\label{n(N)1}
\bn\ci\Psi (N) := \int_0^\infty N(t)\Psi(N(t))\,dt\le C\|w\|\ci{L^\Phi(I)}\,.
\end{align}
It is clear from the proof that one can build the weights for which the left hand side is really much smaller than the right hand side. So our assumptions of the bump conjecture are slightly weaker (for general  weights) than the classical assumptions.

\def\cprime{$'$}
  \def\lfhook#1{\setbox0=\hbox{#1}{\ooalign{\hidewidth\lower1.5ex\hbox{'}\hidewidth\crcr\unhbox0}}}
\providecommand{\bysame}{\leavevmode\hbox to3em{\hrulefill}\thinspace}
\providecommand{\MR}{\relax\ifhmode\unskip\space\fi MR }
\providecommand{\MRhref}[2]{%
  \href{http://www.ams.org/mathscinet-getitem?mr=#1}{#2}
}
\providecommand{\href}[2]{#2}


\begin{thebibliography}{99}

\bibitem{Buc1}{\sc S. Buckley} {\it Summation conditions on weights},  Michigan Math. J. {\bf 40} (1993), no. 1, 153Ð170.

\bibitem{BR} {\sc O. Beznosova, A. Reznikov}, {\it 
Equivalent definitions of dyadic Muckenhoupt and Reverse H\"older classes in terms of Carleson sequences, weak classes, and comparability of dyadic $L\log L$ and $A_\infty$ constants}, arXiv:1201.0520.



\bibitem{CU-M-P-book}  {\sc D. Cruz-Uribe, J. M. Martell,
    C. P\'erez,} {\em Weights, Extrapolation and the Theory of Rubio
    de Francia}, Operator Theory: Advances and Applications, 215,
  Birkhauser, Basel, (2011).

\bibitem{FKP} {\sc R. Fefferman, C. Kenig, J. Pipher}, {\it The theory of weights and the Dirichlet problem for elliptic equations}, Ann. of Math., {\bf 134} (1991) 65--124.


\bibitem{AL}{\sc A. Lerner} {\it On an estimate of Calder\'on--Zygmund operators by positive dyadic operator},  arXiv: 1202.1860v1 [math.CA] 9 Feb 2012.


    \bibitem{NRTV1} {\sc F. Nazarov, A. Reznikov, S. Treil, A. Volberg}, {\em The sharp bump condition for the two-weight problem for classical
singular integral operator: the Bellman function approach}, preprint, Oct. 2011, 1--4.

 \bibitem{NRTV2} {\sc F. Nazarov, A. Reznikov, S. Treil, A. Volberg}, {\em A solution of the  bump conjecture for all Calder\'on--Zygmund operators: the Bellman function approach}, preprint, 9 Feb. 2012, 1--25, arXiv:1202.2406,
 version 1, Feb. 11, version 2,  March 7, 2012.

 
\end{thebibliography}
\end{document}